\documentclass[12pt]{article}

\usepackage{amssymb}
\usepackage{fullpage}

\newcommand{\RR}{{\mathbb{R}}}

  \newtheorem{theorem}{Theorem}

\def\proof{\par{\scshape Proof}. \ignorespaces}
\def\proof{\par{\scshape Proof}. \ignorespaces}
\def\endproof{{\hfill \vbox{\hrule\hbox{%
   \vrule height1.3ex\hskip1.0ex\vrule}\hrule
  }}\par}

\title{A Note on the Eigenvalues of the Google Matrix}
\author{Lars Eld\'en \\
Department of Mathematics, Link\"oping University
}
\date{December 29, 2003}  

\begin{document}
 
\maketitle

\begin{center}
Report LiTH-MAT-R--04-01
\end{center}

\bigskip

Let $P \in \RR^{n \times n}$ be a column-stochastic matrix, i.e. a matrix
with non-negative elements that satisfies $e^T P = e^T$, where $\RR^{1
\times n} \ni e^T = \pmatrix{1 & 1 & \cdots & 1}$. Define 
$$
A = \alpha P + (1 - \alpha) v e^T,
$$
where $0 < \alpha < 1$, and $v \in \RR^n$ is a vector with
non-negative entries that satisfies $e^T v =
1$. Obviously, $A$ is column-stochastic, $e^T A = e^T$, with positive
elements. 

A matrix of this type occurs in the computation of pagerank for the
Google web search engine \cite{brpa:98,pbmw:98}. The pagerank vector is the
(right) eigenvector of $A$ corresponding to the largest eigenvalue in
magnitude, which  is equal to 1.  The corresponding eigenvector has all
non-negative elements, and it is the only eigenvector with this
property. This can be proved using Perron-Frobenius theory, see
e.g. \cite[Chapter 8]{mey:00}.  
 Due to the huge dimension of the matrix,
probably between 3 and 4 billion (December 2003), the only viable
method for computing this eigenvector is the power method, and
variations of it \cite{khg:03,khmg:03a,khmg:03b}  The rate
of convergence of the power method depends on $\lambda_2$,  the second
largest eigenvalue (in magnitude) of $A$, see e.g. \cite[Chapter
7.3]{govl:96}. In \cite{haka:03a} it was 
shown that  $\lambda_2=\alpha$. This result was strenghened in
\cite{lame:03a,lame:03} to that given in Theorem \ref{theo:second}
below. 

The purpose of the present paper is to give a simple alternative proof
of the theorem. 

\begin{theorem}[\cite{lame:03a,lame:03}]\label{theo:second}
Let $P$ be a column-stochastic matrix with  eigenvalues 
$\{1, \lambda_2,  \lambda_3 \ldots, \lambda_n\}$. Then the  eigenvalues of 
$A = \alpha P + (1 - \alpha) v e^T$, where  $0 < \alpha < 1$ and $v$
is a vector with 
non-negative elements satisfying $e^Tv =1$,  are 
$\{1, \alpha \lambda_2, \alpha \lambda_3, \ldots, \alpha \lambda_n\}$.
\end{theorem}

\proof
Define $\hat e$ to be $e$
normalized to Euclidean length 1, and let $U_1 \in \RR^{n \times (n-1)}$
be such that $U=\pmatrix{\hat e & U_1}$ is orthogonal. Then, since 
$\hat e^T P  = \hat e^T$, 
\begin{eqnarray}
U^T P U  &=& \pmatrix{\hat e^T P  \cr U_1^T P} \pmatrix{ \hat e & 
U_1} = \pmatrix{\hat e^T \cr U_1^T P} \pmatrix{\hat e &  U_1} \nonumber\\
\noalign{\vskip 3mm}
&=& \pmatrix{\hat e^T \hat e & \hat e^T  U_1 \cr U_1^T P \hat e &
U_1^T P^T U_1} 
= \pmatrix{1 &  0 \cr  w  & T},\label{eq:simPt}
\end{eqnarray}
where $w =  U_1^T P \hat e $, and $T=U_1^T P^T U_1$. Since we have
made a similarity transformation, the matrix $T$ has the eigenvalues 
$\lambda_2, \lambda_3, \ldots, \lambda_n$. We further have 
$$
 U^T v  = \pmatrix{1/\sqrt{n} \, e^T v \cr  U_1^T v} =
 \pmatrix{1/\sqrt{n} \cr U_1^T v}.
$$
Therefore,
\begin{eqnarray*}
U^T A U &=& U^T (\alpha P + (1 - \alpha)  v e^T)U =
\alpha \pmatrix{1 &  0 \cr  w  & T} + 
(1-\alpha)  \pmatrix{1/\sqrt{n} \cr U_1^T v} \pmatrix{\sqrt{n}\,
  &  0}  \\
\noalign{\vskip 3mm}
&=& \alpha \pmatrix{1 &  0 \cr  w  & T} + 
(1-\alpha) \pmatrix{1 &  0  \cr  \sqrt{n} \, U_1^T v & 0}  =:
 \pmatrix{1 &  0 \cr  w_1  & \alpha T} .
\end{eqnarray*}
The statement now follows immediately.
\endproof

\medskip

The theorem implies that even if $P$ has a multiple eigenvalue equal
to 1, which is actually the case for the Google matrix, the second
largest eigenvalue in magnitude of $A$ is always equal to $\alpha$.

\small

\bibliographystyle{plain}

\end{document}